\documentclass[10pt]{article}
\usepackage{amsthm}
\usepackage{amsmath}
\usepackage{amsfonts}
\usepackage{cancel}
\usepackage{amssymb}

\usepackage[all, knot]{xy}

\def\Z{\mathbb{Z}}

\def\K{\mathbb{K}}

\newcommand{\cale}{{\cal E}}
\newcommand{\calf}{{\cal F}}

\newcommand{\calV}{{\cal V}}

\newcommand{\br}[1]{[\cdot,\cdot]}

\renewcommand{\deg}{{\rm deg}}

\newcommand{\maps}{\colon}

\newcommand{\hh}{\mathfrak{h}}
\newcommand{\g}{\mathfrak{g}}

\newcommand{\cll}{L}
\newcommand{\caa}{A}

\newcommand{\ad}{\operatorname{ad}}

\newcommand{\Hom}{\operatorname{Hom}}

\newcommand{\str}{\operatorname{str}}

\newcommand{\glnv}{{\mathfrak g \mathfrak l}(V)}
\newcommand{\gl}{\mathfrak{gl}}
\newcommand{\half}{\textstyle{\frac{1}{2}}}

\newcommand{\four}{\textstyle{\frac{1}{4}}}
\newcommand{\six}{\textstyle{\frac{1}{6}}}
\newcommand{\eight}{\textstyle{\frac{1}{8}}}
\newcommand{\tw}{\textstyle{\frac{1}{12}}}
\newcommand{\bleft}{[\![}
\newcommand{\bright}{]\!]}

\def\x{{|x|}}
\def\y{{|y|}}
\def\z{{|z|}}
\def\w{{|w|}}

\def\k{{|k|}}
\def\l{{|l|}}

\newcommand{\al}{\alpha}
\newcommand{\be}{\beta}

\newcommand{\om}{\omega}

\def\op{{\oplus}}

\def\la{\langle}
\def\ra{\rangle}

\newcommand{\ot}{\otimes}

\newcommand{\trr}{}

\def\ve#1{(-1)^{#1}}

\newcommand{\doo}{{}_{\bar{0}}}
\newcommand{\di}{{}_{\bar{1}}}

\newcommand{\pf}{\noindent{\bf Proof.}\ }
\newcommand{\allowpagebreak}
\allowdisplaybreaks

\newtheorem{Theorem}{Theorem}[section]
\newtheorem{Proposition}[Theorem]{Proposition}
\newtheorem{Definition}[Theorem]{Definition}

\newtheorem{Lemma}[Theorem]{Lemma}
\newtheorem{Example}[Theorem]{Example}
\newtheorem{Remark}[Theorem]{Remark}

\date{}

\title{Omni-Lie Superalgebras and Lie 2-superalgebras}

\author{Tao ZHANG,\quad Zhangju LIU\\
Department of Mathematics and LMAM\\ Peking University,
Beijing 100871, China\\
}

\begin{document}

\maketitle

\setcounter{section}{0}

\vskip0.1cm

{\bf Abstract}\quad We introduce the notion of omni-Lie superalgebra
as a super version of the omni-Lie algebra introduced by Weinstein.
This algebraic structure gives a nontrivial example of  Leibniz
superalgebra and  Lie 2-superalgebra. We prove that there is a
one-to-one correspondence between Dirac structures of the omni-Lie
superalgebra  and Lie superalgebra structures on subspaces of a
super vector space.



\section{Introduction}
In \cite{Wei}, Weinstein introduced the notion of omni-Lie algebra,
which can be regarded as the linearization of the Courant bracket.
An omni-Lie algebra associated to a vector space $V$ is the direct
sum space $\gl(V)\oplus V$ together with the nondegenerate symmetric
pairing $\la\cdot,\cdot\ra$ and the skew-symmetric bracket
$\bleft\cdot,\cdot\bright$  given by
\begin{eqnarray*}
\la A+u,B+v\ra=\half(Av+Bu),\end{eqnarray*} and
\begin{eqnarray*}
\bleft A+u,B+v\bright=[A,B]+\half(Av-Bu).
\end{eqnarray*}
The bracket $\bleft\cdot,\cdot\bright$ does not satisfy the Jacobi
identity so that  an omni-Lie algebra is not a Lie algebra.
 An omni-Lie algebra is
actually a Lie 2-algebra since Roytenberg and Weinstein proved that
every Courant algebroid gives rise to a Lie 2-algebra
(\cite{Roytenberg}). Rencently, omni-Lie algebras are studied from
several aspects (\cite{BCG}, \cite{InteCourant}, \cite{UchinoOmni})
and are generalized to omni-Lie algebroids and omni-Lie 2-algebras
in \cite{CL, CLS, SLZ}. The corresponding Dirac structures  are also
studied therein.

In this paper, we introduce the notion of omni-Lie superalgebra,
which is the super analogue of omni-Lie algebra. We also study Dirac
structures of omni-Lie superalgebra  in order to characterize   all
possible Lie superalgebra structures on  a super vector space. We
prove that omni-Lie superalgebra  is a Leibniz superalgebra as well
as  a Lie 2-superalgebra, which is a super version of Lie 2-algebra
or a 2-term $L_\infty$-algebra (\cite{Baez, Hue, LS93}).

The paper is organized as follows. In Section 2, we recall some
basic facts for Lie superalgebras. In Section 3, we define omni-Lie
superalgebra on $\cale=\gl(V)\oplus V$ for a super vector space $V$
and study Dirac structures. In Section 4, we prove an  omni-Lie
superalgebra is a Lie 2-superalgebra.

\section{Lie Superalgebras and Leibniz Superalgebras}

We first recall some facts and definitions about Lie superalgebras, basic reference is Kac \cite{Kac}.
We work on a fixed field $\K$ of characteristic 0.

A super vector space $V$ is a $\Z_2$-graded vector space with a
direct sum decomposition $V=V\doo\oplus V\di$.
An element $x\in V\doo\cup V\di$ is called homogeneous. The degree of a homogeneous element $x\in  V_\al,\ \al\in \Z_2$ is defined by $|x|=\al$.
A morphism between
two super vector spaces, $V$ and $W$,is a grade-preserving linear
map:
 $$  f: V \longrightarrow W, \quad f(V_\al)
\subseteq W_\al, \quad \forall \al\in \Z_2.$$ The direct sum
$V\oplus W$ is graded by $$(V\oplus W)\doo = V\doo\oplus W\doo,
\quad (V\oplus W)\di = V\di\oplus W\di,$$ and the tensor product
$V\otimes W$ is graded by
$$(V\ot W)\doo = (V\doo\ot W\doo) \oplus (V\di\ot W\di), \, ~(V\ot W)\di = (V\doo\ot W\di)\oplus (V\di\ot W\doo).$$

\begin{Definition}\label{def:colorLie}  A Lie superalgebra is a super vector space (i.e. $\Z_2$-graded vector space) $\cll=L\doo\op L\di$ together with a bracket
$[\cdot,\cdot]: \cll  \ot\cll  \to\cll $ satisfies,

(i) graded condition: $[L_\al,L_\be]\subseteq L_{\al+\be}$
$\forall\al,\be\in \Z_2$,

(ii) super skew-symmetry:
\begin{align}
[x,y]+\ve{\x\y}[y, x]=0,
\end{align}

(iii) super Jacobi identity:
\begin{equation}\label{J_1}
J_1 := \ve{\z\x}[[x, y],z]+\ve{\x\y}[[y,z],x]+\ve{\y\z}[[z, x],
y]=0,
\end{equation}
where $x,y,z\in L$ are homogeneous elements of degree $\x,\y,\z$ respectively.
\end{Definition}

 One can rewrite the super Jacobi identity in another form:
\begin{align}\label{J2}
J_2 := -\ve{\z\x}J_1=[x,[y,z]]-[[x,y],z] -\ve{\x\y}[y,[x, z]] = 0,
\end{align}
which will be convenient for our use below.

\begin{Example}\label{ex:colorglnv} Let $\caa=A\doo\op A\di$ be an associative superalgebra with multiplication
$A_\al A_\be\subseteq A_{\al+\be}$ for all $\al,\be\in\Z_2$. Define
the bracket:
\begin{equation}\label{gl}
 [x, y]:=xy- \ve{\x\y}y x,\quad \quad \forall x, y\in A.
\end{equation}
Then $(\caa, [\cdot,\cdot])$
is a Lie superalgebra which is denoted by $A_L$.
\end{Example}

\begin{Example} Let $V=V\doo\op V\di$ be a super vector space,
then we have the general linear Lie superalgebra
$$\gl(V) = \gl(V)\doo \op \gl(V)\di,$$
such that
\begin{align*}
\gl(V)\doo &= \Hom(V\doo,V\doo) \op \Hom(V\di,V\di),\\
\gl(V)\di &= \Hom(V\doo, V\di) \op \Hom(V\di, V\doo),
\end{align*}
and the bracket is given by (\ref{gl}). When $dimV\doo = m$,
$dimV\di = n$, $\glnv$ is usually denoted by $\gl(m|n)$.
\end{Example}

 A homomorphism  between two Lie
superalgebras $(L,[\cdot,\cdot ])$ and $(L',[\cdot,\cdot]')$ is
linear map $\varphi: L \to L'$ such that
$$\varphi(L_\al) \subseteq L'_\al, \,
 \quad  \varphi([x,y]) = [\varphi(x),
\varphi(y)]', \quad \forall x,y\in L, \,  \, \forall \al\in \Z_2.$$
A super vector space  $V = V\doo\op V\di$ is called a module of a
Lie superalgebra $L$ or, equivalently, say that $L$ acts on $V$  if
there is a homomorphism $\rho: L\to \glnv$, i.e.,
\begin{align}\label{colormodule}
\rho([x, y])v = \rho(x)\rho(y)v -\ve{\x\y}\rho(y)\rho(x)v.
\end{align}
For simplicity, one often writes  $xv=\rho(x)v$ to denote such an
action. A new Lie superalgebra  can constructed as follows.

\begin{Proposition}\label{semiproduct}\cite{Sch} Let $L$ be a Lie superalgebra with an action on $V$.
Define a  bracket on $L\oplus V$ by
\begin{align}\label{semidirectproduct}
[x+u,y+v]:=[x,y]_L + x\trr v-\ve{|u|\y} y\trr u.
\end{align}
Then  $(L\oplus V,[\cdot,\cdot])$ becomes   a Lie superalgebra,
denoted by  $L\ltimes V$ and  called semidirect product of  $L$ and
$V$.
\end{Proposition}

In \cite{Loday}, Loday  introduced a new algebraic structure, which
is usually called  Leibniz algebra. Its super version is as follows.
\begin{Definition}\label{def:colorLeibniz}\cite{AAA} A Leibniz superalgebra is a super vector space $\cll=L\doo\op L\di$ together with a morphism
$\circ: \cll  \ot\cll  \to\cll $ satisfying $L_\al\circ L_\be\subseteq L_{\al+\be}$ for all $\al,\be\in \Z_2$,  and the super Leibniz rule:
\begin{align}\label{eqn:Leibniz}
x\circ (y\circ z)=(x\circ y)\circ z+\ve{\x\y}y\circ(x\circ z),
\end{align}
for all homogeneous elements $x,y,z\in L$.
\end{Definition}
By definition, it is easy to see that a Leibniz superalgebra is just
a Lie superalgebra if the operation "$\circ$" is super
skew-symmetric. In this case,  the super Leibniz rule above is
actually  the  super Jacobi identity for $J_2$ given by \eqref{J2}.

\section{Omni-Lie Superlgebras and Dirac Stuctures}

Let $V$ be a super vector space, recall that the space $\glnv\op V$ has a $\Z_2$-grading
$$ \cale:= \glnv\op V=(\gl(V)\doo\op V\doo)\op (\gl(V)\di\op V\di).$$
 Like in \cite{InteCourant},  we define an operation $\circ$ on $\glnv\op V$
 as follows:

\begin{equation}
  \label{eq:Leibnizbracket}
   (A+x)\circ (B+y) = [A,B]+Ay.
\end{equation}
Then we have
\begin{Proposition}
$(\cale,\circ)$ is a Leibniz superalgebra.
\end{Proposition}
\pf To check the super Leibniz rule (\ref{eqn:Leibniz}) holds on
$\cale$ under the operation $\circ$, let $e_1=A+x, e_2=B+y, e_3=C+z$
be  homogenous elements of degree $|A|=|x|$, $|B|=|y|$ and
$|C|=|z|$. By definition, we have
\begin{eqnarray*}
&&\{e_1\circ e_2\} \circ e_3-e_1\circ \{ e_2 \circ e_3\} -\ve{\x\y}e_2\circ \{ e_1 \circ e_3\} \\
&=&([A,B]+Ay)\circ (C+z)-(A+x)\circ ([B,C]+Bz)\\
&&-\ve{\x\y}(B+y)\circ ([A,C]+Az)\\
&=&[[A,B],C]-[A,[B,C]]-\ve{\x\y}[B,[A,C]]\\
&&+[A,B]z-ABz-\ve{\x\y}BAz\\
&=&0,
\end{eqnarray*}
where the equality holds  because $\glnv$ is a Lie superalgebra
acting on $V$. \qed \vspace{3mm}

  Note that the above operation is not super skew-symmetric, we can define a
  super  skew-symmetric bracket on  $\cale =\glnv\oplus V$ as its skew symmetrization:
\begin{eqnarray}\label{eq:bracket}
  \bleft A+x, B+y\bright &\triangleq& \half ( (A+x) \circ (B+y) -  (B+y) \circ
  (A+x))\notag\\
&= &[A,B]+\half\left(Ay-\ve{\x\y}Bx\right),
\end{eqnarray}
and define a $V$-valued  inner product, i.e., a non-degenerated
super symmetric bilinear form:
\begin{equation}  \label{eq:symmetric}
  \langle A+x, B+y\rangle\triangleq\half(Ay+\ve{\x\y}Bx).
\end{equation}

 We call the triple  $(\cale, \bleft\cdot,\cdot\bright,
\la\cdot,\cdot\ra)$ an {\bf omni-Lie  superalgebra}. Without the
factor $1/2$ in   bracket $\bleft \cdot, \cdot \bright$,  this would
be the semidirect product Lie superalgebra for the action of $\glnv$
on $V$ described in Proposition \ref{semiproduct}. With the factor
$1/2$, the bracket does not satisfy the super Jacobi identity, which
leads to the concept of Lie 2-superalgebra defined in the next
section. Next we compute the Jacobiator for this bracket.
\begin{Proposition}\label{prop:homotopy}
For $e_1=A+x, e_2=B+y, e_3=C+z\in \cale$, define
\begin{align*}
T(e_1, e_2, e_3):=&\textstyle{\frac{1}{3}}\{ \ve{\z\x} \langle\bleft e_1,e_2\bright, e_3\rangle+ \ve{\x\y}\langle\bleft e_2,e_3\bright, e_1\rangle\\
&+ \ve{\y\z}\langle\bleft e_3,e_1\bright, e_2\rangle\}.
\end{align*}
Let $J_1$ denote the Jacobiator given in (\ref{J_1}) for the bracket
$\bleft\cdot,\cdot\bright$ on $\cale$, then we have
$$J_1(e_1, e_2, e_3)= T(e_1, e_2, e_3).$$
\end{Proposition}
\pf  We compute  both the sides as follows:
\begin{eqnarray*}
&&J_1(e_1,e_2,e_3)\\
&=&\ve{\z\x}\bleft\bleft A+x, B+y\bright, C+z\bright+\mbox{c.p.}\\
&=&\bleft\ve{\z\x}[A,B]+\half\ve{\z\x}\left(Ay-\ve{\x\y}Bx\right), C+z\bright+\mbox{c.p.}\\
&=&\ve{\z\x}[[A,B],C]+\mbox{c.p.}\\
&&+\half\left(\ve{\z\x}[A,B]z-\half\ve{\z\x}\ve{(\x+\y)\z}C\left(Ay-\ve{\x\y}Bx\right)\right)\\
&&+\half\left(\ve{\x\y}[B,C]x-\half\ve{\x\y}\ve{(\y+\z)\x}A\left(Bz-\ve{\y\z}Cy\right)\right)\\
&&+\half\left(\ve{\y\z}[C,A]y-\half\ve{\z\y}\ve{(\z+\x)\y}B\left(Cx-\ve{\z\x}Az\right)\right)\\
&=&\four\ve{\z\x}ABz-\four\ve{\z\x}\ve{\x\y}BAz+\four\ve{\y\z}CAy\\
&&-\four\ve{\y\z}\ve{\x\y}CBx+\four\ve{\x\y}BCx-\four\ve{\z\x}\ve{\y\z}ACy,
\end{eqnarray*}
\begin{eqnarray*}
  &&T(e_1,e_2,e_3)\\
  &=&\textstyle{\frac{1}{3}}\ve{\z\x}\langle\bleft A+x,B+y\bright, C+z\rangle +\mbox{c.p.}\\
  &=&\textstyle{\frac{1}{3}}\ve{\z\x}\langle[A,B]+\half\left(Ay-\ve{\x\y}Bx\right), C+z\rangle +\mbox{c.p.}\\
  &=&\six\ve{\z\x}\left([A,B]z+\half\ve{(\x+\y)\z}C\left(Ay-\ve{\x\y}Bx\right)\right) +\mbox{c.p.}\\
 &=&\six\ve{\z\x}ABz-\six\ve{\z\x}\ve{\x\y}BAz+\tw \ve{\y\z}CAy\\
&&-\tw \ve{\y\z}\ve{\x\y}CBx+\six\ve{\x\y}BCx-\six\ve{\x\y}\ve{\y\z}CBx\\
&&+\tw \ve{\z\x}ABz-\tw \ve{\z\x}\ve{\y\z}ACy+\six\ve{\y\z}CAy\\
&&-\six\ve{\y\z}\ve{\z\x}ACy+\tw \ve{\x\y}BCx-\tw \ve{\x\y}\ve{\z\x}BAz\\
 &=&\four\ve{\z\x}ABz-\four\ve{\z\x}\ve{\x\y}BAz+\four\ve{\y\z}CAy\\
 &&-\four\ve{\y\z}\ve{\x\y}CBx+\four\ve{\x\y}BCx-\four\ve{\z\x}\ve{\y\z}ACy.
\end{eqnarray*}
Thus, the two sides are equal. \qed\vspace{3mm}

The bracket $\bleft\cdot,\cdot\bright$ does not satisfy the super
Jacobi identity so that  an omni-Lie superalgebra is not a Lie
superalgebra. However,  all possible Lie superalgebra structures on
$V$ can be characterized by means of  the omni-Lie superalgebra.


For a bilinear operation $\om$ on $V$ such that $\om: V_\al\times
V_\be\to V_{\al+\be}$, we define the adjoint operator
$$\ad_\om:V_\al\rightarrow \gl(V)_\al, \quad
 \ad_{\om}(x)(y)=\om(x,y)\in V_{\al+\be}$$ where $x\in V_\al, y\in
V_\be$.
 Then the graph of the adjoint operator:
$$\calf_\om =\{\ad_\om x+x ~; \, \forall x\in V\}\subset \cale = \glnv\op V$$
is  a super subspace of   $\cale$. Denote $\calf_\om^\perp$ the
orthogonal complement of $\calf_\om$ in $\cale$ with respect to the
super symmetric bilinear form $\la\cdot,\cdot\ra$ on $\cale $ given
in \eqref {eq:symmetric}.

\begin{Proposition}
\label{prop:realize}  With the above notations,
$(V,\om)$ is a Lie superalgebra if and only if its graph
$\calf_\om$ is maximal isotropic, i.e. $\calf_\om=\calf_\om^\perp$, and
 is closed with respect to the bracket
$\bleft\cdot,\cdot\bright$.
\end{Proposition}

\pf First we see that
\begin{eqnarray*}
\langle\ad_{\om}(x)+x,
\ad_{\om}(y)+y\rangle&=&\half(\ad_{\om}(x)y+\ve{\x\y}\ad_{\om}(y)x)\\
&=&\half(\om(x,y)+\ve{\x\y}\om(y,x)).
\end{eqnarray*}
 This means that  $\om$ is super skew-symmetric if and only if its graph   is isotropic, i.e.
$\calf_\om\subseteq\calf_\om^\perp$. Moreover, by dimension
analysis, we have $\calf_\om$ is  maximal isotropic.

Next let $[x,y] : =\om(x,y)$, we shall check that the super Jacobi
identity on $V$ is satisfied if and only if $\calf_\om$ is closed
under bracket (\ref{eq:bracket}) on $\cale $. In fact,
\begin{eqnarray*}
\bleft \ad_{\om}(x)+x, \ad_{\om}(x)+y\bright&=&[\ad_{\om}(x),\ad_{\om}(y)]+\half(\ad_{\om}(x)y-\ve{\x\y}\ad_{\om}(y)x)\\
&=&[\ad_{\om}(x),\ad_{\om}(y)]+\half(\om(x,y)-\ve{\x\y}\om(y,x))\\
&=&[\ad_{\om}(x),\ad_{\om}(y)]+\om(x,y).
\end{eqnarray*}
Thus this bracket is closed if and only if
$$[\ad_{\om}(x),\ad_{\om}(y)]=\ad_{\om}(\om(x,y)).$$
In this case, for $\forall z\in V$, we have
\begin{eqnarray*}
&&[\ad_{\om}(x),\ad_{\om}(y)](z)-\ad_{\om}(\om(x,y))(z)\\
&=&\ad_{\om}(x)\ad_{\om}(y)(z)-\ve{\x\y}\ad_{\om}(y)\ad_{\om}(x)(z)-\ad_{\om}(\om(x,y))(z)\\
&=&\ad_{\om}(x)\om(y,z)-\ve{\x\y}\ad_{\om}(y)\om(x,z)-\om(\om(x,y),z)\\
&=&\om(x,\om(y,z))-\ve{\x\y}\om(y,\om(x,z))-\om(\om(x,y),z)\\
&=&[x,[y,z]]-\ve{\x\y}[y,[x,z]]-[[x,y],z]\\
&=&0.
\end{eqnarray*}
This is exactly the super Jacobi identity on $V$. Therefore, the
conclusion follows from Definition \ref{def:colorLie}.
\qed\vspace{3mm}

In \cite{BB09}, quadratic Lie superalgebras are studied for a given
inner product $B$  on $V$. In this case, one has the orthogonal Lie
superalgebra $\mathfrak{o}(V)\subset \gl(V)$ and it is easy to see
that $(V,\om, B)$ is a quadratic Lie superalgebra if and only if
$\omega$ satisfies the two conditions in Proposition
\ref{prop:realize} above as well as $\ad_{\om}x\in \mathfrak{o}(V),
\forall x \in V$.


\begin{Definition}
A Dirac structure $L$ of the  omni-Lie superalgebra $(\gl(V)\oplus
V,\bleft\cdot,\cdot\bright,\la\cdot,\cdot\ra)$ is a maximal
isotropic subspace ($L=L^\perp$) and closed under the bracket
$\bleft\cdot,\cdot\bright$.
\end{Definition}

\begin{Remark}
According to Proposition \ref{prop:homotopy}, for  a Dirac structure
$L$,  we have
$$J_1(e_1,e_2,e_3)=T(e_1,e_2,e_3)=0, \quad  \forall e_i\in L.$$
Thus a Dirac structure is a Lie superalgebra, though  omni-Lie
superalgebra is not for itself. In fact,  a Dirac structure is also
a Leibniz subalgebra under the operation $\circ$.
\end{Remark}

By Proposition \ref{prop:realize}, $(V,\om)$ is a Lie superalgebras
if and only if $\calf_\om$ is a Dirac structure of the omni-Lie
superalgebra $\glnv\op V$. In order to give a general
characterization for all Dirac structures of $\cale$, we adapt the
theory of characteristic pairs developed in \cite{liuDirac} (see
also \cite{SLZ}).

For a maximal isotropic subspace $L \subset \gl(V)\oplus V$, set the
subspace $D=L\cap\gl(V)$. Define $D^0\subset V$ to be the null space
of $D$:
$$
D^0=\{x\in V |\ X(x)=0,\ \forall~X\in D\}.
$$
It is easy to see that $D = (D^0)^0$.
\begin{Lemma}\label{lem:Dirac1} With notations above,
a subspace $L$  is  maximal isotropic if and only if   $L$ is of the
form
\begin{equation}\label{pairL}
L=D\oplus \calf_{\pi|_{D^0}}=\{X+\pi(x)+x~|X\in D,~x\in D^0\},
\end{equation}
where  $\pi:V\to \gl(V)$ is a  super skew-symmetric map.
\end{Lemma}
\pf  In the following, we also denote $\pi(x,y)=\pi(x)(y)\in V$ for
convenience. First suppose that $L$ is given by \eqref{pairL}, then
\begin{eqnarray*}
&&\la X+\pi(x)+x,Y+\pi(y)+y\ra\\
&=&\half\{X(y)+ \pi(x, y)+\ve{\x\y}Y(x)+\ve{\x\y}\pi (y, x)\}\\
&=&\half\{\pi(x,y)+\ve{\x\y}\pi(y,x)\}\\
&=&0,  \quad \quad  \forall X+\pi(x)+x~, \, Y+\pi(y)+y~ \in L,
\end{eqnarray*}
 since  $\pi:V\to \gl(V)$ is  super skew-symmetric  so that $L$ is isotropic. Next we prove that $L$ is maximal
isotropic. For $\forall  Z+z \in L^\perp,$
 $$\la X, Z +z \ra = X(z) =0, ~~  \forall  X\in D ~~\Rightarrow ~~
z\in D^0.$$ Moreover,  $\forall x\in D^0$, the equality belew
\begin{align*}
\la X+\pi(x)+x, C+z\ra&=X(z)+\pi(x)(z)+\ve{\x\z}Cx\\
&=\ve{\x\z}\big(C-\pi(z))(x)=0,
\end{align*}
implies that $C-\pi(z) \triangleq Z\in  D$. Thus
$$C+ z = Z+\pi(z)+z\in
L= D\op \calf_{\pi|_ {D^0}}  ~~ \Rightarrow ~~   L = L^\perp.$$ The
converse part is straightforward  so we omit the details.
\qed\vspace{3mm}

The proof of the following  Lemma is skipped since it is
straightforward  and similar to that in \cite{SLZ}.

\begin{Lemma}\label{lem:Dirac2}
Let $(D, \pi)$ be given above. Then $L$ is a Dirac structure if and
only if the following conditions are satisfied:
\begin{itemize}
\item[\rm(1)]~$ D$ is a subalgebra of $\gl(V)$;
\item[\rm(2)]~$\pi\big(\pi(x, y)\big)-[\pi(x),\pi(y)]\in D$, \quad
$ \forall ~x,y \in D^0$;
\item[\rm(3)]~$\pi(x,y)\in D^0$, \quad
$ \forall ~x,y \in D^0$.
\end{itemize}
\end{Lemma}

Such a pair $(D, \pi)$ is called a { \bf characteristic pair} of a
Dirac structure $L$. By means of the two lemmas above, we can
mention the main result in this section.
\begin{Theorem}
There is a one-to-one correspondence between Dirac structures of the
omni-Lie superalgebra $(\gl(V)\oplus V, \la\cdot,\cdot\ra,
\bleft\cdot,\cdot\bright)$ and Lie superalgebra structures on subspaces of $V$.
\end{Theorem}

\pf For any Dirac structure $L= D\oplus\calf_{\pi|_{ D^0}}$,  a Lie
superalgebra structure on $D^0$  is as follows:
 $$
[x,y]_{ D^0}\triangleq\pi(x, y)\in D^0,\quad \forall~x,y\in
 D^0.
$$
 It easy to see that this is a super skew-symmetric map.
For super Jacobi indentity, we have for all $x, y, z\in  D^0$,
 \begin{eqnarray*}
 [[x,y]_{ D^0},z]_{ D^0}&=&\pi([x,y]_{ D^0})(z)=\pi\big((\pi(x)(y)\big)(z)=[\pi(x),\pi(y)](z)\\
 &=&\pi(x)(\pi(y)(z))-\ve{\x\y}\pi(y)(\pi(x)(z))\\
 &=&[x,[y,z]_{ D^0}]_{ D^0}-\ve{\x\y}[y,[x,z]_{ D^0}]_{ D^0}.
\end{eqnarray*}
Thus we get a Lie superalgebra $( D^0,[\cdot,\cdot]_{ D^0})$.

Conversely, for any Lie superalgebra   $( W, [ \cdot, \cdot ]_W) $
on a subspace $W$ of $V$. Define $D$ by
$$
 D= W^0\triangleq\{X\in\gl( V)|~X(x)=0,\ \ \forall~x\in W\}.
$$ Then  $ D^0=(W^0)^0 = W$.
Since Lie superalgebra structure $[\cdot,\cdot]_ W$ gives  a
 super skew symmetric morphism:
  $$\ad: W\to\gl( W), \quad  \ad_x(y)=[x,y]_ W,$$
we take a super skew symmetric morphism  $\pi: V\rightarrow\gl( V)$,
as an extension of $\ad$ from $W = D^0$ to $V$. Thus we get a
maximal isotropic subspace $L=D\oplus\calf_{\pi|_ W}$  from the pair
$(D, \pi)$ as in Lemma \ref{lem:Dirac1}.

 We shall prove that $L$ is a Dirac
structure. Firstly, $\forall X,Y\in  D$ and $x\in W$, we have
$$
[X,Y](x)=XY(x)-\ve{|X||Y|}Y X(x)=0,
$$
which implies that $ D$ is a subalgebra of $\gl( V)$.

Next step is to  prove that $L$ is closed under the bracket
$\bleft\cdot,\cdot\bright$. Remember that $\pi|_W =\ad$ and
$[\cdot,\cdot]_ W$ satisfies the super Jacobi identity, we obtain
$$[\ad_x,\ad_y]=\ad_{[x,y]_ W}=\ad_{\ad_x y} ,\quad\forall~x,y\in W.$$
For any $X\in  D$ and $x,y\in W$, we have
$$
[X,\ad_x](y)=X([x,y]_ W)-\ve{|X|\x}[x,X(y)]=0,
$$
thus $[X,\ad_x]\in D$. On the other hand, we have
\begin{eqnarray*}
&&\bleft X+\ad_x+x,Y+\ad_y+y\bright\\
&=&[X,Y]+[X,\ad_y]+[\ad_x,Y]+[\ad_x,\ad_y]+\half(\ad_x(y)-\ve{\x\y}\ad_y(x))\\
&=&[X,Y]+[X,\ad_y]+[\ad_x,Y]+\ad_{[x,y]_ W}+[x,y]_ W\\
&\in & D\oplus\calf_{\pi|_ W},
\end{eqnarray*}
Thus, we conclude that $L$ is a Dirac structure. Finally, it is easy
to see  that the  Dirac structure  $L$ is independent of the choice
of  extension $\pi$. This completes the proof. \qed\vspace{3mm}

\section{Lie 2-superalgebras}
The concept of Lie $n$-superalgebras is introduced  in \cite{Hue}.
In particular, the axiom of a Lie 2-superalgebra can be expressed
explicitly as follows:
\begin{Definition} \label{2termliealgebra}
A Lie 2-superalgebra
$\calV=(\calV^1\stackrel{d}{\longrightarrow}\calV^0,l_2,l_3)$
consists of the following data:
\begin{itemize}
  \item two super vector spaces $\calV^{0}$ and
   $\calV^{1}$ together with a morphism
   $d\maps \calV^{1} \rightarrow \calV^{0}$;

  \item a morphism $l_{2}=[\cdot, \cdot]\maps \calV^{i} \ot \calV^{j}
   \rightarrow \calV^{i+j}$;

  \item a morphism $l_{3}\maps \calV^{0} \ot \calV^{0} \ot \calV^{0} \rightarrow \calV^{1};$
\end{itemize}
such that, $\forall$ $x,y,z,w\in \calV^{0}$;  $\forall$ $h, k \in
\calV^{1},$
\begin{itemize}
  \item[(a)] $[x,y] +\ve{\x\y}[y,x]=0$;
  \item[(b)] $[x,h] +\ve{\x|h|}[h,x]=0$;
  \item[(c)] $[h,k]=0$;
  \item[(d)] $l_{3}(x,y,z)$ is totally super skew-symmetric;
  \item[(e)] $d([x,h]) = [x,dh]$;
  \item[(f)] $[dh,k] = [h,dk]$;
  \item[(g)]
  $d(l_{3}(x,y,z))=-[[x,y],z]+[x,[y,z]]+\ve{\y\z}[[x,z],y]$;
  \item[(h)]
  $l_{3}(x,y,dh)=-[[x,y],h]+[x,[y,h]]+\ve{|y||h|}[[x,z],h]$;
\item[(i)] $\delta l_3(x,y,z,w) := [x, l_3(y, z,w)] - \ve{\x\y} [y, l_3(x, z, w)]$
\begin{align*}\label{3cocycle}
&+\ve{(\x+\y)\z}[z, l_3(x,y,w)]-[l_3(x, y, z),w]-l_3([x, y], z,w)\notag\\
& + \ve{\y\z}l_3([x,z],y,w)-\ve{(\y+\z)\w}l_3([x,w],y, z)\notag\\
&- l_3(x,[y, z], w) +\ve{\z\w}l_3(x,[y,w],z)- l_3(x,y,[z,w])=0.
\end{align*}
\end{itemize}
\end{Definition}
This is the super analogue of a 2-term $L_\infty$-algebra which is
equivalent to a Lie 2-algebra (see \cite{Baez} for more details).
Here we use the terminology Lie 2-superalgebra instead of 2-term
$L_\infty$-superalgebra.


Now, for a super vector space $V$, let
$$\calV^0=\glnv\op V, \quad
\calV^1=V, \quad d= i: V\hookrightarrow \glnv\op V,$$ where $i$ is
the inclusion map and define operations:
\begin{align*}
l_2=\bleft\cdot,\cdot\bright,\quad l_3=-(-1)^{\z\x}T.
\end{align*}


\begin{Theorem}\label{omniare2term}
With notations above, the omni-Lie superalgebra  $(\cale,
\bleft\cdot,\cdot\bright, \la\cdot,\cdot\ra)$ defines  a Lie
2-superalgebra $(V\stackrel{d}{\hookrightarrow}\glnv\op V,\, l_2, \,
l_3)$.
\end{Theorem}
\pf  For Condition $(a)$, by the grading in $\glnv\op V$ we have
$\deg(A+x)=\deg (A)=\deg (x)$, then
\begin{eqnarray*}
&&\bleft A+x, B+y\bright+\ve{\x\y}\bleft B+y, A+x\bright\\
&=& [A,B]+\half (Ay-\ve{\x\y}Bx )+\ve{|A||B|}[B,A]\\
&&+\ve{\x\y}\half (Bx-\ve{\y\x}Ay )\\
&=& [A,B]+\ve{\x\y}[B,A]+\half (Ay-\ve{\x\y}Bx )\\
&&+\half (\ve{\x\y}Bx-Ay )\\
&=&0.
\end{eqnarray*}
 Conditions $(b)$, $(c)$, $(e)$ and $(f)$ are easy to be
checked.
By Proposition \ref{prop:homotopy}, we have $l_3=J_2$, thus
Conditions $(g)$--$(h)$ hold. For Condition $(i)$, we first verify a
special case by taking $e_1=A, e_2=B, e_3=C, e_4=w$. In fact, by the
definition of $l_3$ and Proposition \ref{prop:homotopy}, we get
\begin{eqnarray*}
&&[A, l_3(B, C,w)] - \ve{\x\y} [B, l_3(A, C,w)] \\
&&+\ve{(\x+\y)\z}[C, l_3(A, B,w)]-[l_3(A, B, C),w]\\
&& - l_3([A, B], C,w)+ \ve{\y\z}l_3([A, C],B,w)-\ve{(\y+\z)\w}l_3([A,w],B, C)\\
&&- l_3(A,[B, C], w) + \ve{\z\w}l_3(A,[B,w], C)- l_3(A,B,[C,w])\\
&=&-\eight A[B,C]w +\eight\ve{\x\y}B[A,C]w-\eight\ve{(\x+\y)\z}C[A, B]w+0\\
&& + \four[[A, B],C]w- \four\ve{\y\z}[[A, C],B]w+\four \ve{\x\y}\ve{\x\z}[[B,C],A]w \\
&&+\eight\ve{\x(\y+\z)}[B,C]Aw-\eight \ve{\y\z}[A,C]Bw+\eight [A,B]Cw\\
&=&  \four\{[[A,B],C]+\ve{\x(\y+\z)}[[B,C],A]+\ve{(\x+\y)\z}[[C,A],B]\}w \\
&&-\eight\{A[B,C] - \ve{\x\y}B[A,C]+ \ve{(\x+\y)\z}C[A, B]\\
&&- \ve{\x(\y+\z)}[B,C]A+ \ve{\y\z}[A,C]B- [A,B]C\}w\\
&=&  \four\{[[A,B],C]-[A,[B,C]]+\ve{\x\y}[B,[A,C]]\}w \\
&&-\eight\{[A,[B,C]] - \ve{\x\y}[B,[A,C]]- [[A,B],C]\}w\\
&=&0.
\end{eqnarray*}
 The general case  can be
checked similarly. \qed\vspace{3mm}

 A Lie 2-superalgebra is called strict if $l_3=0$. This kind of Lie
 2-superalgebras can be described in terms of crossed module.
\begin{Definition} A crossed module of Lie superalgebras consists of a pair of Lie superalgebras $(\g,[\cdot,\cdot]_{\g})$ and $(\hh,[\cdot,\cdot]_{\hh})$ together with an action
of $\g$ on $\hh$ and a homomorphism  $\varphi: \hh\to \g$ such that
$$
\varphi(x\trr h) = [x, \varphi(h)]_{\g}, \quad \varphi(h)\trr k =
[h, k]_{\mathfrak{h}}, \quad \forall h, k\in\hh, \,  \forall x\in\g.
$$
\end{Definition}

\begin{Proposition}
Strict Lie 2-superalgebras are in one-to-one correspondence with
crossed modules of Lie superalgebras.
\end{Proposition}

\pf Let $\calV^1\stackrel{d}{\longrightarrow} \calV^0$ be a strict
Lie 2-superalgebra. Define $\g=\calV^0$, $\hh=\calV^1$, and the
following two brackets on $\g$ and $\hh$:
\begin{eqnarray*}
&&[h,k]_{\hh}=l_2(dh,k)=[dh,k],\quad\forall~ h, k\in \hh=\calV^1;\\
&&[x,y]_{\g}=l_2(x,y)=[x,y],\quad\forall~x, y\in \g=\calV^0.
\end{eqnarray*}
Obviously, $(\g,[\cdot,\cdot]_{\g})$ is a Lie superalgebra by $(a)$
and $(g)$ in Definition \ref{2termliealgebra}.
 By Condition $(h)$, we have
\begin{eqnarray*}
 && -[[h,k]_\hh,l]_\hh+\ve{\k\l}[[h,l]_\hh k]_\hh+[h,[k,l]_\hh]_\hh\\
 &=&-[d[dh,k],l]+\ve{\k\l}[d[dh,l],k]+[dh,[dk,l]]\\
 &=&-[[dh,dk],l]+\ve{\k\l}[[dh,dl],k]+[dh,[dk,l]]\\
 &=&0.
\end{eqnarray*}
This means that  $( \hh, \,[\cdot,\cdot]_\hh )$  is also Lie
superalgebra.  By Condition (e) and taking  $\varphi=d$, we have
$$
\varphi([h,k]_\hh)=d([dh,k])=[dh,dk]=[\varphi(h),\varphi(k)]_\g,
$$
which implies that $\varphi$ is a homomorphism of Lie superalgebras.
Next we define an action of $\g$ on $\hh$ by
$$x
h \triangleq l_2(x,h)=[x,h]\in\hh,$$ which is an action because the
equality,
\begin{eqnarray*}
&&[x,y]\trr h-x\trr (y\trr h)+\ve{\x\y}y\trr (x\trr h)\\
&=&[[x,y],h]-[x,[y,h]]+\ve{\x\y}[y,[x,h]]\\ &=& 0,
\end{eqnarray*}
holds  by Condition $(h)$. Finally, it is easy to check that
\begin{eqnarray*}
&&\varphi(x\trr h)=d([x,h]) = [x,dh]= [x, \varphi(h)]_{\g}\\
&&\varphi(h)\trr k =[dh,k]=[h, k]_{\hh}.
\end{eqnarray*}
Therefore, we obtain a crossed module of Lie superalgebras.

Conversely,  a crossed module of Lie superalgebras
 gives rise to a Lie 2-superalgebra with $d=\varphi$,
$\calV^0=\g$, $\calV^1=\hh$,  $l_3=0$ and the following operations:
\begin{eqnarray*}
~ l_2(x,y)&\triangleq&[x,y]_{\g},\quad \forall
~x,y\in\g;\\
~l_2(x,h)&\triangleq&x h,\quad\forall~ x\in\hh;\\
~l_2(h,k)&\triangleq&0.
\end{eqnarray*}
All of the conditions for a Lie 2-superalgebra can be verified
directly from the definition of a crossed module. \qed\vspace{3mm}

Another kind of  Lie 2-superalgebras is called { \bf skeletal} if
$d=0$. As pointed in \cite{Hue}, Skeletal Lie
2-superalgebras are in one-to-one correspondence with quadruples
$(\g, V, \rho, l_3)$ where $\g$ is a Lie superalgebra, $V$ is a
super vector space, $\rho$ is a representation of $\g$ on $V$ and
$l_3$ is a 3-cocycle on $\g$ with values in $V$. See
\cite{Sch} for more details of the cohomology of Lie superalgebras.

\begin{Example} {\em
Given a quadratic Lie superalgebra $(\g,\, [\cdot,\cdot],\, B)$ over
$\K$, where $B$ is the supertrace $B(x,y)=\str(xy)$ by Kac
\cite{Kac}, a skeletal Lie 2-superalgebra  can be constructed as
follows: $\calV^1=\K, \, \calV^0=\g, \, d=0$ and
\begin{equation}\label{eqn:l2l3string}
l_2(x,y)=[x,y],\quad l_2(x,h)=0,\quad l_3(x,y,z)=B([x,y],z)
\end{equation}
where $x,y,z\in\g, h\in \K$. Condition $(i)$ is from the fact that
Cartan 3-form $l_3$ is closed. That is, by the invariance of $B$,
we have
\begin{eqnarray*}
&&\delta l_3(w,x,y,z)\\
&=& B(w,[x,[y,z]])+ B(w,[x,[y,z]])+\ve{\y\z}B(w,[[x,z],y]) \\
&&+\ve{(\x+\y)\z}B(w,[z,[x,y]])-B(w, [[x,y],z])- \ve{\x\y}B(w,[y,[x,z]])\\
&=&2B(w,[x,[y,z]]+\ve{\y\z}[[x,z],y]-[[x,y],z])\\
&=&0,
\end{eqnarray*}
which holds from  the  super Jacobi identity on $\g$ . Therefore,
$(\K\stackrel{0}{\longrightarrow}\g,l_2,l_3)$ is a Lie
2-superalgebra as a super version of the
 string Lie algebra. }\end{Example}



\begin{thebibliography}{Ra85}
\bibitem{AAA}
Albeverio S, Ayupov S A, Omirov B A. On nilpotent and simple Leibniz
algebras. Comm in Algebra, 2005, 33: 159-172


\bibitem{Baez}
Baez J, Crans A S. Higher-dimensional algebra VI: Lie
 2-Algebras. Theo Appl Cat, 2004, 12:492--528

\bibitem{BB09} Barreiro E, Benayadi S. Quadratic symplectic Lie superalgebras and Lie bi-superalgebras.
J Algebra, 2009, 321: 582--608

\bibitem{BCG}
Bursztyn H, Cavalcanti G,  Gualteri M. Reduction of Courant algebroids and generalized complex structures.
Adv Math,  2007, 211: 726--765



\bibitem{CL}Chen Z, Liu  Z J. Omni-Lie algebroids. J Geom Phys, 2010, 60:799-808

\bibitem{CLS} Chen Z, Liu Z J, Sheng Y. Dirac structures of omni-Lie algebroids. Int J Math, 2011, 22: 1163--1185

\bibitem{Kac} Kac V G. Lie superalgebras. Adv Math, 1977, 26: 8--96


\bibitem{Hue}  Huerta J. Division Algebras, Supersymmetry and Higher Gauge Theory.
Ph.D. thesis,  University of California, 2011. arXiv:1106.3385


\bibitem{InteCourant}
Kinyon K,  Weinstein A. Leibniz algebras, Courant algebroids, and
multiplications on reductive homogeneous spaces. Amer. J. math.,
2001, 123: 525-550


\bibitem{LS93} Lada T, Stasheff J. Introduction to sh Lie algebras for physicists. Int J Theor Phys, 1993, 32: 1087--1103


\bibitem{liuDirac}
Liu Z J. Some remarks on {D}irac structures and {P}oisson reductions.
In {\em Poisson geometry ({W}arsaw, 1998)}, volume~51 of {\em Banach
  Center Publ.}, pages 165--173. Polish Acad. Sci., Warsaw, 2000.

\bibitem{Loday}
Loday J L.  Une version non commutative des alg\`{e}bres de Lie:
les alg\`{e}bres de Leibniz. Enseign Math, 1993, 39: 269--293


\bibitem{Roytenberg}
Roytenberg D,  Weinstein A. Courant Algebroids and Strongly Homotopy
Lie Algebras. Lett Math Phys, 1998, 46: 81--93




\bibitem{Sch} Scheunert M. The Theory of Lie Superalgebras. Lecture Note in Mathematics 716, Springer, 1979.



\bibitem{SLZ}
Sheng Y, Liu Z J, Zhu C C.  Omni-Lie 2-algebras and their Dirac structures. J Geom Phys, 2011, 61: 560--575

\bibitem{UchinoOmni}
 Uchino K. Courant brackets on noncommutative algebras and omni-Lie
 algebras. Tokyo J Math, 2007, 30: 239--255


\bibitem{Wei} Weinstein A. Omni-Lie Algebras. RIMS K\^{o}ky\^{u}roku  2000, 1176: 95--102



\end{thebibliography}
\end{document}